\documentclass[12pt]{amsart} 
\usepackage{amsmath,amssymb,amscd,amsthm,fullpage, graphicx}
\usepackage{latexsym}
\usepackage[usenames,dvipsnames,svgnames,table]{xcolor}
\usepackage{todonotes}
\reversemarginpar

\newtheorem{theorem}{Theorem}[section]

\newtheorem{prop}[theorem]{Proposition}
\newtheorem{definition}[theorem]{Definition}

\newtheorem{corollary}[theorem]{Corollary}

\newtheorem{lemma}[theorem]{Lemma}

\numberwithin{subcase}{case}

\numberwithin{subsubcase}{subcase}
\newcommand{\beginproof}{\noindent{\bf Proof: }}

\newcommand{\EE}{{\mathcal E}}
\newcommand{\FF}{{\mathcal F}}

\newcommand{\K}{{\mathcal K}}

\newcommand{\card}{\mbox{card}}

\newcommand{\PP}{{\mathcal P}}
\newcommand{\PPP}{{\mathbb P}}
\newcommand{\OO}{{\mathcal O}}

\def\RR{{\mathbb R}}
\def\N{{\mathbb N}}

\def\conv{\hbox{\rm conv}}

\def\Fix{\hbox{\rm Fix}\,}
\def\endproof{\begin{flushright}
$ \Box $ \\
\end{flushright}}

\begin{document}
\title{Polytopes of Maximal Volume Product}

\author{Matthew Alexander}
\thanks{The first author is supported in part by the Chateaubriand Fellowship of the Office for Science \& Technology of the Embassy of France in the United States and The Centre National de la Recherche Scientifique funding visiting research at Universit\'e Paris-Est Marne-la-Vall\'ee}

\author{Matthieu Fradelizi}
\thanks{The Second author is supported in part by the Agence Nationale de la Recherche, project GeMeCoD (ANR 2011 BS01 007 01) and the B\'ezout Labex of Universit\'e Paris-Est}

\author{Artem Zvavitch}
\thanks{The third author is supported in part by the U.S. National Science Foundation Grant DMS-1101636 and the B\'ezout Labex of Universit\'e Paris-Est}

\address{Department of Mathematics, Kent State University,
Kent, OH 44242, USA} \email{malexan5@kent.edu}

\address{Universit\'e Paris-Est Marne-la-Vall\'ee,
Laboratoire d'Analyse et de Math\'{e}matiques Appliqu\'ees (UMR 8050)
Cit\'e Descartes - 5, Bd Descartes, Champs-sur-Marne
77454 Marne-la-Vall\'ee Cedex 2, France}\email{matthieu.fradelizi@u-pem.fr }

\address{Department of Mathematics, Kent State University,
Kent, OH 44242, USA} \email{zvavitch@math.kent.edu}

\title{Polytopes of Maximal Volume Product}

\begin{abstract}
For a convex body $K \subset \RR^n$, let $K^z = \{y\in{\mathbb R}^n : \langle y-z, x-z\rangle\le 1, \mbox{\ for all\ } x\in K\}$ be the polar body of $K$ with respect to the center of polarity $z \in \RR^n$. The goal of this paper is to study the maximum of the volume product
 $\mathcal{P}(K)=\min_{z\in {\rm int}(K)}|K||K^z|$, among convex polytopes $K\subset {\mathbb R}^n$ with a number of vertices bounded by some fixed integer $m \ge n+1$.
In particular, we prove that the supremum is reached at a simplicial polytope with exactly $m$ vertices and we provide a new proof of a result of Meyer and Reisner showing that, in the plane, the regular polygon has maximal volume product among all polygons with at most $m$ vertices. Finally, we  treat the case of polytopes with $n+2$ vertices in $\RR^n$.
\end{abstract}
\maketitle

\section{Introduction}

We denote the inner product of two vectors $x, y \in \RR^n$ by $\langle x, y \rangle$  and the length of a vector $x \in \RR^n$ by $|x|$.
A {\it convex body} is a compact convex subset of $\RR^n$ with non empty interior.  We say that a set $K$ is {\it symmetric} if it is centrally symmetric with center at the origin, i.e. $K=-K$.
We write $|A|$ for the $k$-dimensional Lebesgue measure (volume)  of a measurable set $A \subset \RR^n$, where $k$ is the dimension of the minimal affine subspace containing $A$. We denote by $\conv(A)$ the {\em closed} convex hull of a set $A \subset \RR^n$, by ${\rm int}(A)$ its interior and by $\conv(A,B, C, \dots)$ the closed convex hull of $A\cup B\cup C, \dots$. For $a,b\in \RR^n$, we denote $[a,b]$ the segment joining $a$ to $b$: $[a,b]= \{ (1-t)a+tb: t \in [0,1]\}$.  We will frequently refer to \cite{Gr}, \cite{Sc} and \cite{Z} for general references for convex bodies and polytopes and their properties.

 The polar body $K^z$ of $K$ with the center of polarity $z$ is defined by
$$
K^z = \{y\in\RR^n : \langle y-z, x-z\rangle\le 1 \mbox{\ for all\ } x\in K\}
.$$
If the center of polarity is taken to be the
origin, we denote the polar body of $K$ by $K^\circ$ .
Note that $K^z=(K-z)^\circ + z$, and the bipolar theorem says that $(K^z)^z =K$, for $z\in {\rm int}(K)$ (see \cite{Gr}, p. 47).

 A well known
result of Santal\'o \cite{S} (see also  \cite{Sc}, p. 546) states that in every convex body $K$ in
$\RR^n$, there exists a unique point $s(K)$, called the
\emph{Santal\'o point} of $K$, such that
$$
|K^{s(K)}| = \min_{z\in {\rm int}(K)} |K^z|.
$$

The \emph{volume product} of $K$ is defined by
$$
\PP(K) = \inf \{|K| |K^z| : z\in {\rm int}(K)\}= |K|\ |K^{s(K)}|.
$$

The volume product is affinely invariant, that is,
$\PP(A(K))=\PP(K)$ for every affine isomorphism $A: \RR^n \rightarrow
\RR^n$. Observe that if we denote $L=K^{s(K)}$ then
$$\PP(K^{s(K)})=|L||L^{s(L)}|\le |L||L^{s(K)}|=|K^{s(K)}||K|=\PP(K).$$
The set of all convex bodies in $\RR^n$ is compact with respect to the Banach-Mazur distance and $K\mapsto \PP(K)$ is continuous in Hausdorff distance (see, for example, \cite{FMZ}), so that it is natural to ask for maximal and minimal values of $\PP(K)$.
The Blaschke-Santal\'o inequality states that
$$
\PP(K) \le \PP(B^n_2),
$$
where $B^n_2$ is the Euclidean unit ball. The  equality in the above inequality is possible only for ellipsoids (\cite{S}, \cite{P}, see \cite{MP}
or also \cite{MR2} for a simple proof of both the inequality and
the case of equality).

The minimality of $\PP(K)$ is an open question, often called Mahler's conjecture \cite{Ma1,Ma2}, which states that, for every convex body $K$ in $\RR^n$,
\begin{equation}\label{eq:inverse_santalo}
\PP(K) \ge \PP(\Delta^n)=\frac{(n+1)^{n+1}}{(n!)^2},
\end{equation}
where $\Delta^n$ is an $n$-dimensional simplex. It is also conjectured
that equality in (\ref{eq:inverse_santalo}) is attained only if
$K$ is a simplex.
The  symmetric case of Mahler conjecture states that for every {\it symmetric} convex body $K \subset \RR^n$:
\begin{equation}\label{eq:inverse_santalo_sym}
\PP(K) \ge \PP(B_1^n)=\PP(B_\infty^n)=\frac{4^{n}}{n!},
\end{equation}
where $B_1^n$ and $B_\infty^n$ are the  cross-polytope and its dual, the cube, respectively.

The inequalities (\ref{eq:inverse_santalo})  and (\ref{eq:inverse_santalo_sym}) for
$n=2$ were proved by Mahler \cite{Ma1} with the case of equality
proved by Meyer \cite{Me2} in the general case and by Reisner \cite{R1} in the symmetric case. Other cases, such as bodies of
revolution, were treated in \cite{MR1}. Several special cases of the conjecture, most of them for symmetric bodies, can be found in
\cite{BF, BMMR, SR, R1, GMR, Me1, R2, FMZ, NPRZ, Ki, KiR, RSW, GM}.
A special case of   $n$ dimensional polytopes with at most $n+ 3$ vertices (or facets) was treated  in \cite{MR2}. The proof of this last result is based on the method of shadow systems which we shall elaborate on, applying it to finding the maximum of the volume product.

Observe that  an isomorphic version of reverse Santal\'o inequality  was proved by Bourgain and Milman
\cite{BM}:
$$\PP(K) \ge c^n\PP(B^n_2),$$
where $c$ is a positive constant; Kuperberg \cite{Ku}  gave a new proof of this
result with a better constant (see also \cite{Na}, \cite{GPV} for different proofs of the inequality and \cite{AGM}, \cite{RZ} for more information).

The goal of this paper is to study the maximal value of the volume product when we restrict ourselves to the class of polytopes with a bounded number of vertices. We start by introducing the basic tools in section \ref{sec:shadows}. We recall the definition and properties of shadow systems introduced by Rogers and Shephard and state the propositions of Campi, Gronchi, Meyer, and Reisner connecting shadow systems to the volume product. In Theorem \ref{sup-attained} we show that the maximum value of the volume product among all convex polytopes in $\RR^n$ with $m$ vertices is increasing in $m$. Next, in Theorem \ref{sup-simplicial} we prove that the polytopes of maximal volume product among polytopes with at most $m$ vertices must satisfy some identities which imply in particular that it is simplicial.

 In Section~\ref{sec:polygons} we give a new proof of the result of Meyer and Reisner \cite{MR3} showing that the regular $m$-gon is the only $m$-gon with maximal volume product among polygons with at most $m$ vertices. 

Then, in section \ref{sec:n+2}, we consider the case of convex polytopes with $n+2$ vertices in $\RR^n$ and in Theorem \ref{th:n+2 points} we prove that the polytope with maximal volume product is the convex hull of two simplices living in supplementary affine subspaces of dimensions $\lceil \frac{n}{2} \rceil$ and $\lfloor \frac{n}{2} \rfloor$.


\section{Main tools}\label{sec:shadows}

The main tool in the proof of our results is the technique of shadow systems  of convex sets introduced by Rogers and Shephard \cite{RS} and generalized by Shephard \cite{Sh} in the following way. Let $C$ be a closed convex set in $\RR^{n+1}$. Let $(e_1, \ldots, e_{n+1})$ be an orthonormal basis of $\RR^{n+1}$, we write $\RR^{n+1}=\RR^n\oplus\RR e_{n+1}$, so that $\RR^n=e_{n+1}^\bot$. For every $x\in\RR^n$ let $P_x$ be the projection onto $\RR^n$ parallel to $e_{n+1}-x$: for $z\in \RR^n$ and $s\in \RR$,
$$P_x(z+se_{n+1})= z + sx.$$ 
We denote $K_x=P_x(C)\subset\RR^n$. Let $I$ be a convex subset of $\RR^n$. Then we say that the family $(K_x)_{x\in I}$ is a shadow system of convex sets. Rogers and Shephard \cite{RS} proved that $x\mapsto|K_x|$ is convex on $I$. Campi and Gronchi \cite{CG} proved that if moreover the convex bodies $K_x$ are origin symmetric, for every $x\in I$, then $x\mapsto|K_x^\circ|^{-1}$ is convex on $I$.
In \cite{MR2}, Meyer and Reisner generalized this result to the non-symmetric case and studied the equality case. The following proposition is one of our key tool:

\begin{prop}[\cite{MR2}]\label{lm:mr}
Let $I$ be a convex subset of $\RR^n$ and $(K_x)_{x\in I}$, be a  shadow system of convex bodies in $\RR^n$ then $x\mapsto \left|K_x^{s(K_x)}\right|^{-1}$ is convex on $I$. 
\end{prop}

As a corollary, if the volume of $K_x$ is constant,  then $x\mapsto\PP(K_x)^{-1}$ is convex. Moreover if the function $x\mapsto |K_x|$ is affine then $x\mapsto \PP (K_x)$ is the quotient of an affine function by a convex one. As it was noticed in \cite{MR1} Lemma 12 and in  \cite{FMZ} Corollary 2, it follows that it is quasi-concave: i.e.  $\{x\in I: \PP (K_x)\ge s\}$ is convex, for every $s>0$. 

The following (classical) lemma is a useful observation for us to treat the maximal cases of the volume product:
\begin{lemma}\label{lem:max:conv}
Let $K \subset \RR^n$ be a convex body and $F: K \to \RR$ be a concave continuous function. Assume that $K$ and $F$ are  invariant under linear isometries  $T_1,...,T_m$. Then there is $x_0 \in K$ such that $T_i(x_0)=x_0$, for all $i=1,\dots, m$ and $F(x_0) \ge F(x)$ for all $x \in K$.
\end{lemma}
\beginproof Let us first assume that the function $F$ is strictly concave, i.e. $F((x+y)/2) > (F(x)+F(y))/2$, for $x \not =y$.   Then by continuity of $F$ and compactness of $K$, the maximum of $F$ is reached at $x_0 \in K$, moreover this point is unique by strict concavity, indeed if $x,y \in K$ are two distinct maximums, then $F((x+y)/2) > F(x)$ and $(x+y)/2 \in K$. 

Moreover the function $F$ is invariant under a map $T_i$ so then $F(T_ix_0)=F(x_0)$. But because the maximum is reached at unique point we have $T_ix_0=x_0$. 

Now if $F$ is only concave and not necessary strictly concave, 
we may approximate $F$ by a sequence of strictly concave functions $F_k(x)=F(x)-|x|^2/k$.
The maps $T_i$ are isometries and thus $F_k(T_ix)=F_k(x)$ for all $i\in 1, \dots, m$ and $k \in \mathbb{N}$.
By the previous argument applied to $F_k$, we deduce that for each $k$ there is a unique $x_k\in K$ 
such that $\max_{x\in K} F_k(x)=F_k(x_k)$ and $T_ix_k=x_k$ for all $i \in 1,\dots,m$. 
Now since $K$ is compact we may select a convergent  subsequence $\{x_{k_l}\}$ of $\{x_k\}$. Let  $\lim x_{k_l}=x_0$, then $x_0\in K$ and 
by continuity of $T_i$,we get $T_ix_0=x_0$ for all $i$. Moreover, by continuity of $F$ we get $\lim F(x_{n_l}) =  F(x_0)$, therefore $\max_{x\in K} F(x)=F(x_0)$.
\endproof

In a number of  places throughout the paper we will often say that two bodies are ``close enough'', this is measured with respect to the Hausdorff distance: for two non-empty subsets $K, L \subset \RR^n$ we define their Hausdorff distance  $d_H(K, L)$ by 
$$
d_{H}(K,L) = \max\{\,\sup_{x \in K}  \inf_{y \in L} |x-y|,\, \sup_{y \in L}  \inf_{x \in K}|x-y|\,\}.
$$

Finally, the following proposition is a combination of Propositions 1 and 2 of Kim and Reisner \cite{KiR} which will help us to estimate the behaviour of $|L^z|$ when $z$ is close enough to $s(L)$. 
 \begin{prop}[\cite{KiR}]\label{prop:kr}
Let $K$ and $L$ be two convex bodies in $\RR^n$. Then there exists $\delta(K)$ such that, 
if $d_H(K,L)\le \delta(K)$ then $$|L^{s(L)}|=|L^{s(K)}|+O(d_H(K,L)^2),$$
where $O$ depends only on $K$.
As a consequence, if the Santal\'o point of $K$ is at the origin and $d_H(K,L)\le \delta(K)$ then 
$$\PP(L)=|L||L^\circ|+O(d_H(K,L)^2).$$
\end{prop}

\section{Properties of Polytopes of Maximal Volume Product}\label{char}

\begin{definition}
For $n\ge 1$  we denote by ${\mathcal K}^n$ the set of all convex bodies in $\RR^n$ endowed with the Hausdorff distance. For $m\ge n+1$, we denote by $\PPP^n_m$ the  subset of $\mathcal K^n$ consisting of the polytopes in $\RR^n$ with non empty interior  
having at most $m$ vertices and by $\PPP^n=\cup_{m\in\N}\PPP^n_m$, the dense subset of $\mathcal K^n$ consisting of all polytopes with non-empty interior. We denote by $M^n_m$ the supremum of the volume product of polytopes with at most $m$ vertices and non empty interior in $\RR^n$
$$M^n_m:=\sup_{K\in\PPP^n_m}\mathcal P(K).$$
\end{definition}

Recall that from Blaschke-Santal\'o inequality one has $\sup_{K\in\K^n}\PP(K)=\PP(B_2^n)$. By the continuity of the function $K\mapsto \mathcal P(K)$  on $\mathcal K_n$ (see for example Lemma 3 in \cite{FMZ}) and the density of $\PPP^n$ in $\K^n$ we deduce that $\lim_{m\to+\infty}M^n_m=\PP(B_2^n)$. Our aim is now to establish that the sequence $M^n_m$ is strictly increasing. We start with a lemma that is of independent interest and gives a better understanding on the behavior of the  volume product functional. 

\begin{lemma}\label{Kx}
Let $n,m\in\N$ with $m\ge n+1$ and $K\in \PPP^n_m$. Let $F$ be a facet of $K$ with exterior normal $u\in S^{n-1}$, let $x_F$ be in the relative interior of $F$ and let $K_t=\conv(K, x_F+tu)$, for $t>0$. Then for $t$ small enough the volume product of $K_t$ is strictly larger than the volume product of $K$:
$$\PP(K_t)>\PP(K).$$
\end{lemma}

Notice that the polytope $K_t$ defined in the above proposition has exactly $m+1$ vertices.

\beginproof 
We may assume that the Santal\'o point of $K$ is at the origin. Let $h>0$ such that the affine hyperplane spanned by $F$ is $H=\{x:  \langle x,  u\rangle =h\}$ and $K\subset H^-$, where $H^{-}=\{x : \langle x,  u \rangle \le h\}$. Let $F_1,\dots, F_k$ be the facets of $K$ which are adjacent to $F$ and for $1\le i\le k$, denote by $u_i$ the exterior normal of $F_i$. Let  $h_i>0$ be such that $H_i=\{x : \langle x, u_i\rangle =h_i\}$ is the spanned affine hyperplane of $F_i$. Thus 
$$
K\subset \bigcap\limits_{1\le i\le k}H_i^-.
$$
We also denote by 
$$R=\{x : \langle x,  u \rangle \ge h,  \ \langle x,  u_i \rangle \le h_i,\ \forall 1\le i\le k\}$$ 
the polyhedral region bounded by $F$ and the $H_i$, $i=1, \dots, k$. For every $x\in R$, let $K_x=\conv(K,x)$ then $(K_x)_{x\in R}$ is a shadow system and 
 $$|K_x|=|K|+\frac{1}{n}|F|(\langle x,  u\rangle-h).$$
Notice that $x\mapsto |K_x|$ is affine and thus $(K_x)_{x\in R}$ is an affine volume shadow system. 
 Let $x_F$ be an interior point of $F$ and let  $x_t=x_F+tu$, then if $t>0$ and small enough we get $x_t\in R$. Moreover, using that  $x_F\in F$ and thus $\langle x_F,  u\rangle =h$ we get $\langle x_t, u\rangle =h+t$ and
 $$|K_{x_t}|=|K|+\frac{t}{n}|F|.$$
 By polarity, the point $u/h$ is a vertex of $K^\circ$, the points $u_i/h_i$ are its adjacent vertices and $K_x^\circ=\{y\in K^\circ; \langle y, x \rangle \le 1\}$ is the truncation of the polytope $K^\circ$ by the halfspace $\{y: \langle y,  x\rangle \le 1\}$. For every $x$ in the interior of $R$ this truncation cuts off the vertex $u/h$ of $K^\circ$. It also cuts the edges $[u/h ; u_i/h_i]$ at some points $v_i=(1-\lambda_i) u/h+\lambda_i u_i/h_i$, where $\lambda_i\in[0,1]$ is determined by the fact that $\langle v_i,  x\rangle =1$. This gives 
 $$\lambda_i=\frac{(\langle x,  u\rangle -h)h_i}{\langle x, (h_i u-h u_i)\rangle}.$$
 Thus 
 $$v_i-\frac{u}{h}=-\lambda_i\left(\frac{u}{h}-\frac{u_i}{h_i}\right)=-\frac{\langle x, u\rangle -h}{h}\times \frac{h_i u-h u_i}{\langle x, h_i u-h u_i\rangle }. $$
 Moreover one has 
 $$K^\circ\setminus K_x^\circ=\conv\left(\frac{u}{h},v_1,\dots, v_k\right)
 =\frac{u}{h}+\conv\left(0,v_1-\frac{u}{h},\dots, v_k-\frac{u}{h}\right).$$
 Hence 
\begin{eqnarray*}
|K_x^\circ|&=&|K^\circ|-\left|\conv\left(0,v_1-\frac{u}{h},\dots, v_k-\frac{u}{h}\right)\right|\\
&=&|K^\circ|-\left(\frac{\langle x, u\rangle -h}{h}\right)^n\left|\conv\left(0,\frac{h_1 u-h u_1}{\langle x,  (h_1 u-h u_1)\rangle },\dots, \frac{h_k u-h u_k}{\langle x,  (h_k u-h u_k)\rangle }\right)\right|.
\end{eqnarray*}
Applying this for $x=x_t$ and using that $\langle x_t, u\rangle =h+t$, we get
$$
|K_{x_t}^\circ|=|K^\circ|-\left(\frac{t}{h}\right)^n\left|\conv\left(0,\frac{h_1 u-h u_1}{\langle x_t, (h_1 u-h u_1)\rangle},\dots, \frac{h_k u-h u_k}{\langle x_t, (h_k u-h u_k)\rangle }\right)\right|.
$$
Thus for $t$ small enough, we obtain 
$$|K_{x_t}^\circ|=|K^\circ|+O(t^n).$$
Hence 
$$|K_{x_t}||K_{x_t}^\circ|=(|K|+t|F|/n)\left(|K^\circ|+O(t^n)\right)
=|K||K^\circ|+t|K^\circ||F|/n+o(t).$$
Moreover, one has $d_H(K,K_{x_t})\le c(K) t$ for some constant $c(K)$ depending on $K$ only. Hence it follows from Proposition \ref{prop:kr} that for $t>0$ small enough 
$$\PP(K_{x_t})=|K_{x_t}||K_{x_t}^\circ|+O(t^2).$$
We deduce that 
$$\PP(K_{x_t})= \PP(K)+t|K^\circ||F|/n+o(t)>\PP(K).$$
\endproof

\noindent{\bf Remark:} It is tempting to state Lemma \ref{Kx} in a stronger form, saying that for any $n$-dimensional polytope $K \subset \RR^n$ and a point $x \in \RR^n$, such that $\conv(K, x)$ has more vertices than $K$ one has $\PP(\conv(K,x)) \ge \PP(K)$. But such a statement is wrong. This can be seen by a direct computation, or from the following observation: consider $K=B^2_\infty$ and $x_\epsilon=(10, 1-\epsilon)$. Then, the continuity of the volume product gives us
$$
\lim\limits_{\epsilon \to 0} \PP(\conv\{B^2_\infty, x_\epsilon\})=\PP(\conv\{
(1,-1); (-1,-1); (-1, 1); (10,1)\}) < \PP(B^2_\infty),
$$
where the last inequality follows from direct computation (see also Theorem \ref{th:polygons}, below).

\begin{theorem}\label{sup-attained}
Let $n\ge 1$ and $m\ge n+1$. The supremum $M^n_m$ is achieved at some polytope with exactly $m$ vertices and the sequence $M^n_m$ is strictly increasing in $m$.
\end{theorem}

\beginproof The fact that the supremum $M^n_m$ is achieved follows the proof of the corresponding statement on the infimum established, for example, in Proposition 2 and Lemma 4 of \cite{FMZ}.  
 By the affine invariance of $\mathcal P$ and F. John's theorem (see \cite{Sc}, page 588) one has 
$$M^n_m:=\sup_{K\in\PPP^n_m}\mathcal P(K)=\sup\{\mathcal P(K): K\in\PPP^n_m, B_2^n\subset K\subset nB_2^n\}.$$
Note that   $\{K\in\PPP^n_m: B_2^n\subset K\subset nB_2^n\}$ is compact in Hausdorff distance. Moreover the function $K\mapsto \mathcal P(K)$ is continuous on $\mathcal K_n$ (see for example Lemma 3 in \cite{FMZ}). Therefore as the supremum of a continuous function $\mathcal P$ on a compact set, we conclude that the supremum $M^n_m=\sup_{K\in\PPP^n_m}\mathcal P(K)$ is attained at some polytope $K_m$ with at most $m$ vertices. 

Now let us prove that any polytope $K_m$ achieving the supremum has exactly $m$ vertices. The proof goes by induction on $m$. For $m=n+1$, the result is clear. Let $m\ge n+1$ be fixed and assume that the result is known for $K_m$. So $K_m$ has exactly $m$ vertices. From Lemma \ref{Kx} there exists $x$ outside $K_m$ such that $K_m(x)=\conv(K_m,x)$ has a volume product strictly larger than $K$.
Since $K_m(x)\in \PPP^n_{m+1}$, it follows that 
$$M^n_{m+1}=\PP(K_{m+1})\ge \PP(K_m(x))>\PP(K_m)=M^n_m.$$
We conclude that $K_{m+1}$ has exactly $m+1$ vertices and that the sequence $m\mapsto M^n_m$ is strictly increasing.
\endproof

\noindent{\bf Remark:} Notice that since the Euclidean ball is known to be the maximum in volume product among all bodies, then from this and the above theorem we can see that there is no polytope which is a local maximum of the volume product among all convex bodies.

Recall that one says that a polytope is simplicial if all its facets are simplices.

\begin{theorem}\label{sup-simplicial}
Let $n\ge 1$ and $m\ge n+1$. Let $K$ be of maximal volume product among polytopes with at most $m$ vertices. Then $K$ is a simplicial polytope.
\end{theorem}

For the proof, we need to introduce some more notation concerning polytopes. For any polytope $K$ we denote by $\EE(K)$ the set of its vertices and by $\FF(K)$ the set of its facets.\\

\beginproof 
Let $K$ be a polytope with the origin in its interior. For any facet $F\in\FF(K)$, we denote $u_F$ its exterior normal and by $h_F$ its distance to the origin. Let $x$ be a vertex of $K$. Denote by $\FF(x)$ the set of facets of $K$ containing $x$.  We denote by $F_x$ the facet of $K^\circ$ corresponding to $x$: 
$$F_x=\{y\in K^\circ; \langle y,x\rangle=1\}=\left\{y\in K^\circ; \left\langle y,\frac{x}{|x|}\right\rangle=\frac{1}{|x|}\right\}.$$
Notice that $F_x$ has $\frac{x}{|x|}$ as its exterior normal and its distance to the origin is $1/|x|$. Now we introduce a modification of $K$ that was used by Meyer and Reisner in \cite{MR3} in the plane: we define $K_t=\conv(K, (1+t)x)$, for small values of $t>0$, so we extend $K$ in the direction of $x$. Then,
$$
|K_t|=|K|+\sum_{F\in\FF(x)}|\conv(F, (1+t)x)|.
$$
For any $F\in\FF(x)$, one has $\langle u_F,x\rangle=h_F$, thus
$$
|\conv(F, (1+t)x)|=\frac{1}{n}|F|(\langle u_F, (1+t)x\rangle-h_F)= \frac{t}{n}|F|h_F=t|\conv(F,0)|.
$$
Hence 
$$
|K_t|=|K|+t\sum_{F\in\FF(x)}|\conv(F,0)|.
$$
The result of this change of $K$ is a cutting for $K^\circ$ parallel to the facet $F_x$:
$$
K_t^\circ=\left\{y\in K^\circ ; \langle y,x\rangle \le \frac{1}{1+t}\right\}= \left\{y\in K^\circ ; \left\langle y,\frac{x}{|x|}\right\rangle \le \frac{1}{(1+t)|x|}\right\}.
$$
For sufficiently small $t>0$ the distance between the facet $F_x$ and the new parallel facet is 
$$d_x=\frac{1}{|x|}\left(1-\frac{1}{1+t}\right)=\frac{t}{(1+t)|x|}=\frac{t}{|x|}+o(t).$$
Thus it is not difficult to see that we get
$$
|K_t^\circ|=|K^\circ|-t\frac{|F_x|}{|x|}+o(t)=|K^\circ|-nt|\conv(F_x,0)|+o(t).
$$
Together, we get
$$
|K_t||K_t^\circ|=|K||K^\circ|+t\left(|K^\circ|\sum_{F\in\FF(x)}|\conv(F,0)|-n|K||\conv(F_x,0)|\right)+o(t).
$$
Now we assume that the Santal\'o point of $K$ is at the origin. Then using Proposition \ref{prop:kr} similarly to Lemma \ref{Kx}, since $d_H(K,K_t)=O(t)$ we get $\PP(K_t)=|K_t||K_t^\circ|+O(t^2)$. Thus, for $t>0$,
\begin{equation}\label{eq:PKt}
\PP(K_t)=\PP(K)+t\left(|K^\circ|\sum_{F\in\FF(x)}|\conv(F,0)|-n|K||\conv(F_x,0)|\right)+o(t).
\end{equation}
Now let us assume that $K$ has maximal volume product among polytopes with at most $m$ vertices.
Since $K_t$ has also $m$ vertices, it follows that $\PP(K_t)\le\PP(K)$ and thus using \eqref{eq:PKt} for any vertex $x$ of $K$ we have
\begin{eqnarray}\label{deriv-class}
|K^\circ|\sum_{F\in\FF(x)}|\conv(F,0)|\le n|K||\conv(F_x,0)|.
\end{eqnarray}
Summing on all the vertices of $K$ we get
$$
\sum_{x\in\EE(K)}|K^\circ|\sum_{F\in\FF(x)}|\conv(F,0)|\le \sum_{x\in\EE(K)}n|K||\conv(F_x,0)|=n|K||K^\circ|.
$$
Simplifying by $|K^\circ|$ and inverting sums in the left hand side gives 
$$
\sum_{F\in\FF(K)}\card(\EE(F))|\conv(F,0)|\le n|K|.
$$
Since for any facet $F$, one has $\card(\EE(F))\ge n$, we get
$$
n|K|\le \sum_{F\in\FF(K)}\card(\EE(F))|\conv(F,0)|\le n|K|.
$$
Thus we get equality in all previous inequalities, which implies that for any facet $F$ one has $\card(\EE(F))= n$. Therefore every facet $F$ is a simplex and so $K$ is simplicial.
We also get the following consequence, for any vertex $x\in\EE(K)$ one has 
\begin{eqnarray}\label{eq:simplicial}
|K^\circ|\sum_{F\in\FF(x)}|\conv(F,0)|= n|K||\conv(F_x,0)|.
\end{eqnarray}
\endproof
\noindent{\bf Remark:}
 Let us notice that if a polytope $K$ minimizes the volume product among polytopes with at most $m$ vertices then the same argument shows that the inequality (\ref{deriv-class}) is reversed: for every vertex $x$ of $K$ one has
$$|K^\circ|\sum_{F\in\FF(x)}|\conv(F,0)|\ge n|K||\conv(F_x,0)|.$$
It's easy to see that simplices and $B_1^n$ satisfy the above inequality. 

One may also establish the following lemma generalizing equation (\ref{eq:simplicial}).
\begin{lemma}\label{lem:eq}
Let $n\ge 1$ and $m\ge n+1$. Let $K$ be of maximal volume product among polytopes with at most $m$ vertices. Assume that the Santal\'o point of $K$ is at the origin. Let $x\in\EE(K)$ be a vertex of $K$ and denote by $\FF(x)$ the facets of $K$ containing $x$. Then one has 
\begin{eqnarray}\label{eq:vector}
|K^\circ|\sum_{F\in\FF(x)}|\conv(F,0)|y_F= n|K||\conv(F_x,0)|g_{F_x},
\end{eqnarray}
where $g_{F_x}$ denotes the center of gravity of the facet $F_x$ of $K^\circ$ corresponding to $x$ and for every $F\in\FF(x)$, 
$y_F$ denotes the vertex of $K^\circ$ corresponding to $F$.
\end{lemma}
\beginproof
From Theorem \ref{sup-simplicial}, we know that $K$ is simplicial. Using Lemma 5 of \cite{FMZ}, we may apply a more general shadow system than the one used in the proof of Theorem \ref{sup-simplicial}. Let $Q=\conv(\EE(K)\setminus\{x\})$ and for $z$ in a neighborhood of $x$ define $K(z)=\conv(Q,z)$. Then one has 
$$
|K(z)|=|K|+\frac{1}{n}\sum_{F\in\FF(x)}|F| \langle z-x,u_F\rangle.
$$
Hence, using that $y_F=u_F/h_F$, 
$$
\nabla |K(z)|_{z=x}=\frac{1}{n}\sum_{F\in\FF(x)}|F| u_F=\sum_{F\in\FF(x)}|\conv(F,0)|y_F.
$$
Next we notice that $K(z)^\circ=\{y\in Q^\circ; \langle y, z\rangle \le 1\}$ and, using formula (3) on page 347 of \cite{MR}, we get  
$$
\nabla |K(z)^\circ|_{z=x}=-\frac{|F_x|}{|x|}g_{F_x}=-n|\conv(F_x,0)|g_{F_x}.
$$
Because all facets of $K$ are simplices, $z$ can move freely in a neighborhood of $x$ and thus for $K$ maximizing the volume product, we get that 
$\nabla \PP(K(z))_{z=x}=0$. Again, from Proposition \ref{prop:kr} one has $\PP(K(z))=|K(z)||K(z)^\circ|+O(|z-x|^2)$ thus 
$$
\nabla (|K(z)||K(z)^\circ|)_{z=x}=\nabla \PP(K(z))_{z=x}=0=|K|\nabla |K(z)^\circ|_{z=x}+|K^\circ|\nabla |K(z)|_{z=x}.
$$
Hence we get that for every vertex $x\in\EE(K)$ 
$$
|K^\circ|\sum_{F\in\FF(x)}|\conv(F,0)|y_F= n|K||\conv(F_x,0)|g_{F_x}.
$$
\endproof

\noindent{\bf Remark:}\\ 
1) Notice that if $K$ has maximal volume product among symmetric polytopes with at most $m$ vertices, then in the proof of Theorem \ref{sup-simplicial} one can consider $K_t=\conv(K, \pm(1+t)x)$ and we get that $K$ satisfies the inequality (\ref{deriv-class}) and thus $K$ must be simplicial.\\   
2) We should note that a simple and simplicial polytope is either a polygon or a simplex (see, for example, \cite{Z}, page 67). Thus if $K$ has maximal volume product among the polytopes with at most $m>d+1$ vertices in dimension $d>2$ then its polar is not of maximal volume product in its class and doesn't necessarily satisfy equation (\ref{eq:vector}).   
Still, following \cite{MR3} we may claim that, in $\RR^2$,  $K^\circ$  will satisfy the combinatorial properties of (\ref{eq:vector}).

Indeed, let $K \subset \RR^2$ be of maximal volume product among polygons with at most $m$ vertices and Santal\'o point at the origin. 
Let $L=K^\circ$, $y$ be a vertex of $L$ and define $L(z)$ in the same way we defined  $K(z)$ in the proof of Lemma
\ref{lem:eq}, i.e. $z$ is a small perturbation of the vertex $y$. Using that $K$ is a polygon we get that  $(L(z))^\circ$ has the same number of vertices as $K$.

We get that $\PP((L(z))^\circ)$ is maximal when $z=y$ and that $\nabla \PP((L(z))^\circ)|_{z=y}=0$.
Now, again, as in the proof of Lemma
\ref{lem:eq} we  use Proposition \ref{prop:kr}, and since $$d_H((L(z))^\circ,K)=d_H(((L(z))^\circ)^\circ,K^\circ)=d_H(L(z),L)=O(|z-y|)$$ we have 
$$
|((L(z))^\circ)^{S((L(z))^\circ)}|=|((L(z))^\circ)^{S(K)}|+O(|z-x|^2)=|L(z)|+O(|z-y|^2).
$$ 
So then $$|L(z)||L(z)^\circ| =|((L(z))^\circ)^{S((L(z))^\circ)}||L(z)^\circ| + O(|z-x|^2)=\PP((L(z))^\circ)+O(|z-x|^2).$$
Thus $\nabla(|L(z)||L(z)^\circ|)_{z=y}=0$ and we can conclude similarly as in proof of Lemma \ref{lem:eq}.

\section{Maximality in $\RR^2$}

Let us fix some notation. For $\theta\in[0,2\pi]$, we set $R_\theta$ to be the rotation about the origin of angle $\theta$ in the oriented plane $\RR^2$. We denote by $e_1, e_2$ the canonical basis of $\RR^2$.
For $m\ge 3$ we consider the regular polygon with $m$ vertices and unit circumcircle:
$$P_m:=\mathrm{conv}\left\{R_\frac{2k\pi}{m}(e_1);\; k=0,\ldots,m-1\right\}.$$
A simple calculation shows that $|P_m|=m\sin(\pi/m)\cos(\pi/m)$. Note that 
$$P_m^\circ=\frac{1}{\cos(\pi/m)}R_\frac{\pi}{m}(P_m)$$
is also a regular polytope (obtained by rotating and dilating $P_m$). We deduce that $|P_m^\circ|=m\tan(\pi/m)$ and the volume product of $P_m$ is thus 
\begin{equation}\label{eq:vrp}
\mathcal P(P_m)=\big( m\sin(\pi/m)\big)^2.
\end{equation}
Notice that $m\mapsto\mathcal P(P_m)$ is an increasing sequence. Indeed, the function $x\mapsto \sin(x)/x$ is positive and decreasing on $[0,\pi)$.

We shall give a new proof of the following result of Meyer and Reisner  \cite{MR3}.

\label{sec:polygons}
\begin{theorem} \label{th:polygons}
Let $m\ge3$ and let $K$ be a polygon in $\RR^2$ with at most $m$ vertices, then 
$$\mathcal P(K)\le \mathcal P(P_m),$$
with equality if and only if $K$ is an affine image of $P_m$.
\end{theorem}

 We start the proof of Theorem \ref{th:polygons} with a lemma showing that if  a polygon $K$ achieves the maximum of the volume product among polygons with a fixed number of vertices, then each vertex of $K$ must be on the line passing through the Santal\'o point of $K$ and the middle of its two adjacent vertices.
\begin{lemma} \label{lem:3vert}
Let $K \in \PPP^2_m$ have maximal volume product among polytopes in $\PPP^2_m$ and Santal\'o point at the origin. Then for any vertex $x$ of $K$ there exists a real number $\lambda=\lambda(x)$ such that $x=\lambda (x_1+x_2)$, where $x_1$ and $x_2$ are the vertices of $K$ adjacent to $x$.
\end{lemma}

\beginproof
Let $x$ be a vertex of $K$ and denote by $x_1$ and $x_2$ its two adjacent vertices. Denote by   $y_1$ and $y_2$ the vertices of $K^\circ$ corresponding to the edges $[x,x_1]$ and $[x,x_2]$ of $K$. We apply equation (\ref{eq:vector}) of Lemma \ref{lem:eq} to our situation, the center of gravity of the edge $F_x^*$ of $K^\circ$ corresponding to the $x$ is the center of gravity of the edge $[y_1, y_2]$, hence it is the middle of the segment $[y_1, y_2]$, thus   $g_{F_x^*}=(y_1+y_2)/2$. So equation (\ref{eq:vector}) becomes: 
$$
|K| |\conv(0,y_1,y_2)|(y_1+y_2)=|K^\circ| (|\conv(0,x,x_1)|y_1+|\conv(0,x,x_2)|y_2).
$$ 
Because these quantities are equal and $y_1$ and $y_2$ are linearly independent, we may identify and get 
$$
|K| |\conv(0,y_1,y_2)|=|K^\circ| |\conv(0,x,x_1)|=|K^\circ| |\conv(0,x,x_2)|.
$$
Choosing an orientation, we deduce that 
$$\det(x_1,x)=|\conv(0,x,x_1)|=|\conv(0,x,x_2)|=\det(x,x_2)=-\det(x_2,x).
$$
Thus $\det(x_1+x_2,x)=0$. Hence 
there exists a real number $\lambda=\lambda(x)$ such that $x=\lambda (x_1+x_2)$.
\endproof

Hence we proved that for a polygon with maximal volume product and Santal\'o point at the origin, each vertex must be a multiple of the sum of its two adjacent vertices. By the second remark after Lemma \ref{lem:eq} we can also conclude that this property holds in the dual as well.  Now we will show that for any polygon with $m$ vertices that has the property of Lemma \ref{lem:3vert} the  constant $\lambda$ is independent of the triple of vertices:

\begin{lemma} \label{lem:constant}
Let $K$ be a convex polygon with $m$ vertices and maximal volume product with its Santal\'o point at the origin.  Then there exists a real number, $\lambda>1/2$ such that for any  vertices, $v_1$, $v_2$, and $v_3$, with $v_2$ adjacent to both $v_1$ and $v_3$ one has $v_2=\lambda (v_1+v_3)$.
\end{lemma}

\beginproof
Let us order the vertices of the polygon counterclockwise as $x_1,\ldots, x_m$ and the vertices of the dual $y_1, \ldots, y_m$ with $y_i$ such that $\langle x_i, y_i\rangle = \langle  x_{i+1}, y_i\rangle=1$. 
By Lemma \ref{lem:3vert} applied to $P$, there exists real numbers $\lambda_i$  so that for all $1\le i\le m-1$
$$
x_i=\lambda_i(x_{i-1}+x_{i+1}).
$$  
Taking the scalar product with $y_i$ and $y_{i-1}$, we get
$$
\langle x_{i-1},y_i\rangle= \langle x_{i+1},y_{i-1}\rangle=\frac{1}{\lambda_i}-1.
$$

Now we can use the Remark after  Lemma \ref{lem:eq} to claim that $P^\circ$ will also satisfy the combinatorial conditions of Lemma \ref{lem:3vert}. Thus, there exists $\mu_i$ such that for all $1\le i\le N-1$
$$
y_i=\mu_i(y_{i-1}+y_{i+1}).
$$
Taking the scalar product with $x_i$ and $x_{i+1}$, we get
$$
\langle x_i,y_{i+1}\rangle= \langle x_{i+1},y_{i-1}\rangle=\frac{1}{\mu_i}-1.
$$
Using the equations above, we deduce that $\lambda_i=\mu_i=\lambda_{i+1}$. 

\endproof
 
Now using Lemma \ref{lem:constant} and standard techniques to solve recurrence relations we can prove Theorem \ref{th:polygons}.

\noindent{\bf Proof of Theorem \ref{th:polygons}:} 
 By an affine transform, we may assume that the Santal\'o point of $P$ is at the origin. Denote by $v_1,\dots, v_m$ the vertices of $P$ (counting clockwise). Again, applying linear transformation we may assume  $v_0=v_m=e_1$, where, as before,  $e_1$ is the first vector of the canonical basis $(e_1,e_2)$ of $\RR^2$. From Lemma \ref{lem:constant} we have the recurrence relation for the vertices: $tv_k=v_{k+1}+v_{k-1}$.  Then the recurrence holds also for the coordinates $x_k$ and $y_k$ of $v_k$. Since $v_k$ is a vertex of $P$, one has $0<t<2$ thus the roots of the equation $y^2-ty+1=0$ are $\alpha=e^{i\theta}$ and $\beta=e^{-i\theta}$, with $\cos(\theta)=t/2$.
Thus there exists $A,B\in\RR$ such that for every $k$ 
$$
y_k=A\cos(k\theta)+B\sin(k\theta),
$$
with initial conditions $y_0=y_m=0$. Since $y_0=0$, we get $A=0$. Notice that if $B=0$ then all $y_k$ are 0 and thus all vertices lie on the $x$-axis, so we discard this possibility. 
 So by the initial conditions we have $B\sin(m\theta)=0$ hence $\sin(m\theta)=0$. Thus there exist $j \in \N$ such that $\theta=\frac{j\pi}{N}$. The first coordinate $x_k$ of $v_k$  satisfies the same recurrence relation so there exists $C$ and $D$ such that for every $k$ 
$$
x_k=C\cos\left(\frac{jk\pi}{m}\right)+D\sin\left(\frac{jk\pi}{m}\right),
$$
with the initial conditions $x_0=x_m=1$. Since $x_0=1$, we get $C=1$. Since $x_m=1$, we get that $j$ must be even.  Taking in to account that $K$ has exactly $N$ vertices we get that $j=2$. Finally
$$
(x_k, y_k)=\left(\cos\left(\frac{2k\pi}{m}\right)+D\sin\left(\frac{2k\pi}{m}\right), B\sin\left(\frac{2k\pi}{m}\right)\right),
$$
and thus $K$ is the linear image of the regular $m$-gon $P_m$ by the map 
$$T= \left(\begin{array}{cc}
1 & D    \\
0 & B   
\end{array}\right).$$ 
\endproof

 \section{Convex hull of $n+2$ points in $\RR^n$} \label{sec:n+2} 

For $K\subset \RR^{n}$ being a convex body, we 
define $\OO(K)=\{T\in O_{n}(\RR) ; TK=K\}$ and $\Fix(K)=\{x\in\RR^n  ; Tx=x, \forall T\in\OO(K)\}$. We shall consider convex bodies $K$ such that $\Fix(K)$ is one point, the origin. In this case, notice that all affine invariant points attached to $K$ coincide with this point. In particular the Santal\'o point of $K$ satisfies $s(K)=0$. 

\begin{theorem}\label{th:sym}
Let $1\le k\le n-1$ be integers and let $E$ and $F$ be two supplementary subspaces in $\RR^n$ of dimensions $k$ and $n-k$ respectively. Let $L\subset E$ and $M\subset F$ be convex bodies of the appropriate dimensions such that $\Fix(L)=\Fix(M)=\{0\}$. Then for every $x\in L$ and  $y\in M$
$$\PP(\conv(L-x,M-y))\le\PP(\conv(L,M))=\frac{\PP(L)\PP(M)}{\binom{n}{k}},$$
with equality if and only if $x=y=0$.
\end{theorem}

\beginproof  Using the invariance of the volume product under linear transformations we may assume that $E$ and $F$ are perpendicular. 
Now consider the following shadow system $(x,y)\mapsto K_{x,y}= \conv(L-x, M-y)$, for $(x,y)\in K\times L$. Computing the volume, we get  
$$
|K_{x,y}|=\frac{|L-x||M-y|}{\binom{n}{k}}=\frac{|L||M|}{\binom{n}{k}}.
$$
 So $(x,y)\mapsto K_{x,y}$ is a volume constant shadow system. Thus, from Proposition \ref{lm:mr}, the function $(x,y)\mapsto \PP(K_{x,y})^{-1}$ is convex on $L\times M$. Moreover, for any $(S,T)\in \OO(L)\times\OO(M)$ one has $(S\times T)(L\times M)=L\times M$ and for any $(x,y)\in L\times M$ 
$$
K_{S(x),T(y)}=\conv(L-S(x),M-T(y))=\conv(S(L-x), T(M-y))=(S\times T)(K_{x,y}).
$$
Thus $\PP(K_{S(x),T(y)})=\PP(K_{x,y})$.
This means that the function $(x,y)\mapsto \PP(K_{x,y})$ is invariant under the action of $\OO(L)\times\OO(M)$. By Lemma \ref{lem:max:conv}, we deduce that its maximum occurs at a fixed point of $\OO(L)\times\OO(M)$, which is reduced to the origin by the hypotheses. The equality case is clear.
\endproof

\begin{corollary}\label{cor:sym}
Let $L\subset \RR^{n-1}$ be a convex body such that $\Fix(L)$ is one point. Then among all double pyramids $K=\conv(L, x, y)$ in $\RR^n$ with base $L$ separating apexes $x$ and $y$, the volume product $\PP(K)$ is maximal when $x$ and $y$ are symmetric  with respect to the Santal\'o point of $L$.
\end{corollary}



\begin{theorem} \label{th:n+2 points}
Let $K$ be the convex hull of $n+2$ points. Let $q=\lfloor\frac{n}{2}\rfloor$ and $p=\lceil\frac{n}{2}\rceil=n-q$.
Then 
$$
\PP(K)\le \frac{(p+1)^{p+1}(q+1)^{q+1}}{n!p!q!},
$$
with equality if and only if $K$ is the convex hull of two simplices $\Delta_{q}$ and $\Delta_{p}$ living in supplementary affine subspaces of dimensions $q$ and $p$ respectively.
\end{theorem}

\beginproof
Let $K$ be a body in $\PPP^n_{n+2}$.  Then by Radon's theorem there exists $1\le k\le n-1$ such that one may split the $n+2$ vertices of $K$ into two subset $I$ and $J$, with $\card(I)=k+1$ and $\card(J)=n+1-k$ in such a way that if $L=\conv(I)$ and $M=\conv(J)$ then $L\cap M\neq\emptyset$. Since $K=\conv(I,J)=\conv(L,M)$ is full dimensional in $\RR^n$, it follows that $L$ and $M$ are non-degenerate simplices in supplementary affine subspaces $E$ and $F$ of dimension $k$ and $n-k$. By affine invariance, we may assume that $E\cap F=\{0\}$, $E$ and $F$ are orthogonal to each other and $L$ and $M$ are standard simplices of their respective dimensions so that one may write $L=\Delta_k+s(L)$, where $s(L)$ is the Santal\'o point of $L$ and $\Delta_k$ is a regular simplex of dimension $k$ with Santal\'o point at the origin. In the same way, one has $M=\Delta_{n-k}+s(M)$. Since $\Fix(\Delta_k)=\Fix(\Delta_{n-k})=0$, we may apply Theorem \ref{th:sym} to $\Delta_k$ and $\Delta_{n-k}$. We get that 
$$
\PP(K)\le\PP\left(\conv(\Delta_k,\Delta_{n-k})\right)=\frac{\PP(\Delta_k)\PP(\Delta_{n-k})}{\binom{n}{k}}:=f_n(k).
$$
Recall that the volume product of a non-degenerate simplex $\Delta_n$ in $\RR^n$ is 
$$\PP(\Delta_n)=\frac{(n+1)^{n+1}}{(n!)^2}.$$
Hence after simplification, we get 
$$
f_n(k)= \frac{1}{n!}\times\frac{(k+1)^{k+1}}{k!}\times\frac{(n-k+1)^{n-k+1}}{(n-k)!}=\frac{g(k)g(n-k)}{n!},
$$
where $g(x)=\frac{(x+1)^{x+1}}{\Gamma(x+1)}$, for $x\ge 0$. Then with the change of variable $t=xu$ we get, for $x>0$, 
$$\frac{\Gamma(x)}{x^x}=\frac{1}{x^x}\int_0^{+\infty}e^{-t}t^{x-1}dt=\int_0^{+\infty}(ue^{-u})^x\frac{du}{u}.
$$
Hence from H\"older's inequality the function $x\mapsto \frac{\Gamma(x)}{x^x}$ is log-convex on $(0,+\infty)$. It follows that $g$ is log-concave on $\RR_+$. So $f_n$ is log-concave as well, and since it satisfies $f_n(k)=f_n(n-k)$, for all $0\le k\le n$, we deduce that 
$$
f_n(k)=\frac{g(k)g(n-k)}{n!}\le \frac{g(\lfloor\frac{n}{2}\rfloor)g(\lceil\frac{n}{2}\rceil)}{n!},
$$
with equality if and only if $K=\conv(\Delta_{\lfloor\frac{n}{2}\rfloor},\Delta_{\lceil\frac{n}{2}\rceil})$.
\endproof

{\bf Remark:} One may conjecture that for $k\le n$, among polytopes with at most $n+k$ vertices, the convex hull of $k$ simplices living in orthogonal subspaces of dimensions $\lfloor\frac{n}{k}\rfloor$ or $\lceil\frac{n}{k}\rceil$ have maximal volume product (see Gluskin-Litvak \cite{GL} where such bodies where considered). Theorem \ref{th:n+2 points} establishes this conjecture for $k=2$.

\end{document}